\documentclass[a4paper,11pt]{article}
\usepackage{indentfirst,latexsym,bm,color}
\usepackage{amsmath,amssymb,amsfonts}

\usepackage[top=1in,left=1in,right=1in]{geometry}

\makeatletter

\@addtoreset{equation}{section} \makeatother

\begin{document}
\newtheorem{dingli}{Theorem}[section]
\newtheorem{dingyi}[dingli]{Definition}
\newtheorem{tuilun}[dingli]{Corollary}
\newtheorem{zhuyi}[dingli]{Remark}
\newtheorem{yinli}[dingli]{Lemma}

\title{ The minus order for idempotents
\thanks{ This work was supported by NSF of
China (Nos: 11671242, 11571211,11601339) and the
Fundamental Research Funds for the Central Universities
(GK201801011).}}
\author{  \ Yuan Li$^a$\thanks{E-mail address:
 liyuan0401@aliyun.com},\ \ \ Jiajia Niu$^a$,\ \ \ Xiaoming Xu$^b$\thanks{E-mail address:
 xuxiaoming2620@aliyun.com}}
\date{} \maketitle\begin{center}
\begin{minipage}{16cm}
{ \small $$    \ \ \ a.   \ School \ of \ Mathematics \ and \
Information \ Science,\ Shaanxi \ Normal \ University, $$ $$ Xi'an,
710062, People's \ Republic \ of\ China. $$ }{\small $$ \ \ \ b.  \ School \ of \ Science,\ Shanghai \ Institute \ of \ Technology,  $$ $$ Shanghai, \ 20418, \ People's \
Republic \ of \ China. $$}
\end{minipage}
\end{center}
 \vspace{0.05cm}
\begin{center}
\begin{minipage}{16cm}
\ {\small {\bf Abstract }
Let $P$ and $Q$ be idempotents on a Hilbert space $\mathcal{H}.$
The minus order $P\preceq Q$ is defined by the equation $PQ=QP=P.$ In this note,
 we first  present some necessary
and sufficient conditions for which the supremum and infimum of idempotents $P$ and $Q$ exist with respect
to the minus order.
Also, some properties of the minimum $Q^{or}$ are characterized, where
$Q^{or}$=min $\{P^{'}: P^{'}$ is an orthogonal projection on $\mathcal{H}$ with $Q \preceq P^{'} \}.$
\\}
\endabstract
\end{minipage}\vspace{0.10cm}
\begin{minipage}{16cm}
\end{minipage}
\end{center}
\begin{center} \vspace{0.01cm}
\end{center}

\section{\textbf{Introduction and preliminaries}}
Let $\mathcal{H}$ and $\mathcal{K}$ be separable complex Hilbert spaces,
 $\mathcal{B(H,K)}$
be the set of all bounded linear operators from
$\mathcal{H}$ into $\mathcal{K}.$
An operator $A\in\mathcal{B(H)}$ is called positive, if
$A\geqslant 0,$ that is $\langle
Ax,x\rangle \geqslant0 $
 for all $x\in\mathcal{H},$ where $ \langle , \rangle$ is the inner product of $\mathcal{H}.$
 Also, we denote by $\mathcal{B({H})}^{+}$ the set of all positive bounded linear operators on $\mathcal{H}.$ For
  $A \in\mathcal{B({H})}^{+},$ $A^{\frac{1}{2}}$ is the positive square root of $A.$ In particular, $|A|: =(A^{*}A)^{\frac{1}{2}}$ is the absolute value of operator $A,$
where $A^{*}$ is the adjoint operator of $A.$

For an operator $T\in \mathcal{B(H,K)},$  $N(T),R(T)$ and $\overline{R(T)}$ denote the null
 space, the range of $T,$ and the closure of $R(T),$ respectively.
Let $\mathcal{M}$ and $\mathcal{N}$  be closed subspaces of $\mathcal{H}.$ We write by $\mathcal{M}+ \mathcal{N}$ the linear subspace spanned by $\mathcal{M}$ and $\mathcal{N}.$ When $\mathcal{M}\cap\mathcal{N}=\{0\},$
we denote by $\mathcal{M}\dotplus\mathcal{N}=\mathcal{M}+ \mathcal{N}$ the direct sum of $\mathcal{M}$ and $\mathcal{N}.$ In Particular, $\mathcal{M}\oplus\mathcal{N}$ is the orthogonal sum and $\mathcal{M}\ominus\mathcal{N}
=\mathcal{M}\cap(\mathcal{M}\cap\mathcal{N})^\perp$
is the orthogonal minus.
  Also, $P_{\mathcal{M}}$ denotes the orthogonal projection onto
the closed subspace $\mathcal{M}$ and a rank-one operator $x\otimes y$ is defined by $(x\otimes y)z=\langle z,y\rangle x$
for all $z\in\mathcal{H}.$
Moreover, $\overline{\{x\}}$ represents the one-dimensioned subspace spanned by a nonzero vector $x\in\mathcal{H}$
and $\overline{\{x_1,x_2\}}=\overline{\{x_1\}}\vee\overline{\{x_2\}}$ for nonzero vectors $x_1,x_2\in\mathcal{H}.$

An operator $J\in \mathcal{B(H)}$ is said to be a symmetry (or self-adjoint unitary operator) if $J=J^*=J^{-1}.$
In this case, $J^+=\frac{I+J}{2}$ and $J^-=\frac{I-J}{2}$ are mutually annihilating orthogonal projections. If
$J$ is a non-scalar symmetry, then
 an indefinite inner product is defined by \[[x,y]:=\langle Jx, y\rangle \qquad(x,y\in \mathcal{H})\] and
 $(\mathcal{H}, J)$ is called a Krein space [1].
Let $\mathcal{B(H)}^{Id}$ and $\mathcal{P(H)}$ be the set of all idempotents and orthogonal projections on $\mathcal{H},$ respectively.
For $P\in\mathcal{B(H)}^{Id},$ if
$ran(P) = M$ and $ker(P) = N,$ then $P$ is called the idempotent operator onto $M$ along $N.$
An idempotent $P\in\mathcal{B(H)}^{Id}$ is said to be a $J$-projection,
 if $P=JP^{*}J.$ The
existence of $J$-(positive) projections and its properties are studied in [12-15].

As usual, the operator order (Loewner partial order) $A\leq B$ between two bounded self-adjoint operators
is defined as $A-B\leq 0.$
For $P, Q \in \mathcal{B(H)}^{Id},$ we write $P\preceq{Q}$ if $PQ=QP=P.$  This relation defines a partial order on $\mathcal{B(H)}^{Id}.$ Indeed, it follows from [5, Definition 3.1 or 16 Definition 1] that
this partial order is equivalent to the minus partial order which is confined to $\mathcal{B(H)}^{Id}.$
The minus partial order is a well known order defined and studied for matrices and later on for operators acting on Hilbert spaces by many authors (see [3,5,16,17]).
It is trivial that $P\preceq Q$ if and only if $P\leq Q$ for $P,Q \in \mathcal{P(H)}.$
For $P,Q \in \mathcal{B(H)}^{Id},$ we denote by $P\underset{\tiny{\preceq}}{\Large{\vee}} Q$ the supremum,  equivalently, the least upper bound of
$P$ and $Q$ with respect to the partial order $\preceq,$ if it exists. To be more precise,
 $P\underset{\tiny{\preceq}}{\Large{\vee}} Q$ is an idempotent, uniquely determined
by the following properties: $P\preceq P\underset{\tiny{\preceq}}{\Large{\vee}} Q,$ $Q\preceq P\underset{\tiny{\preceq}}{\Large{\vee}} Q$ and if
$Q'\in \mathcal{B(H)}^{Id}$ satisfies both $P\preceq Q'$ and $Q\preceq Q',$ then $P\underset{\tiny{\preceq}}{\Large{\vee}} Q \preceq Q'.$
Analogously, $P\underset{\tiny{\preceq}}{\Large{\wedge}} Q$
denotes the greatest lower bound of $P$ and $Q$ with respect to the order $\preceq.$

Let $Q \in \mathcal{B(H)}^{Id}.$ In the following Proposition 3.3, we show that the sets of
$\{P: P\preceq Q \hbox{ and } P\in \mathcal{P(H)}\}$
and $\{P: Q\preceq P \hbox{ and } P\in \mathcal{P(H)}\}$
have the maximum and minimum with respect to the order $\preceq,$ respectively. Denote by
$$ Q_{or}:=\underset{\tiny{\preceq}}{\Large{max}} \{P: P\preceq Q \hbox{ and } P\in \mathcal{P(H)}\}$$
and
$$ Q^{or}:=\underset{\tiny{\preceq}}{\Large{min}} \{P: Q\preceq P \hbox{ and } P\in \mathcal{P(H)}\}.$$

Suppose that $\{Q_{n}\}_{n\in \mathbb{N}}$ is a sequence in $\mathcal{B(H)}.$
$\{Q_{n}\}_{n\in \mathbb{N}}$ is said to converge in the WOT topology to $Q$
(denote by $Q_{n}\xrightarrow[n\rightarrow \infty]{WOT}Q$ )
if $\langle Q_{n}x,y\rangle\xrightarrow[n\rightarrow \infty] \ \langle Qx,y\rangle$
for every $x,y\in \mathcal{H}.$
For $\{Q_{n}\}_{n\in \mathbb{N}} \in \mathcal{B(H)}^{Id},$  we denote
 $Q_{n}{\overset{\tiny{WOT}}{\Large{\nearrow}}}Q,$ if $Q_{n}\xrightarrow[n\rightarrow \infty]{WOT}Q$
 and $Q_{n}\preceq Q_{n+1},$ for all $n=1,2,\cdots.$
 Analogously, we write $Q_{n}{\overset{\tiny{WOT}}{\Large{\searrow}}}Q$ if $Q_{n}\xrightarrow[n\rightarrow \infty]{WOT}Q$
 and $Q_{n+1}\preceq Q_{n}.$

For a given partial order of $\mathcal{B(H)},$
studying its lattice properties is an interesting problem.
That is equivalent to giving the necessary and sufficient conditions for the existence of
supremum and infimum for two arbitrary operators with respect to this partial order.
For the operator order and the star partial order,
the existence of infimum and supremum have been studied in different contexts
 (see Refs. [2,6-11,19]). However,
 for the minus partial order of $\mathcal{B(H)},$
 studying its lattice properties seems difficult.
The conditions for which
the supremum and infimum of $\mathcal{B(H)}$ with respect to the minus partial order exist
 have not yet been discovered. In this note, we shall make some attempts
  in this topic. In Section 2, we mainly consider the lattice properties of $\mathcal{B(H)}^{Id}$ with respect to the minus partial order.
We present the necessary
and sufficient conditions for which $P\underset{\tiny{\preceq}}{\Large{\vee}} Q $ exists and characterize the specific structures of
$P \underset{\tiny{\preceq}}{\Large{\vee}} Q$ if it exists.
In Section 3, we first give the existence of $Q_{or}$ and $ Q^{or}.$ Then we extend a similar result for the $J$-projections. That is, we get that $Q_{or}$ and $ Q^{or}$ are $J$-projections, if $Q$ is a $J$-projection.
Also, if $P\in \mathcal{P(H)}$ is a $J$-projection, we present the equivalent condition under which there is a $J$-projection $Q\in \mathcal{B(H)}^{Id} \backslash \mathcal\{{\mathcal{P(H)}}\}$
such that $Q^{or}=P.$

\section{ Conditions for the existence of $P\underset{\tiny{\preceq}}{\Large{\vee}} Q$ and $P\underset{\tiny{\preceq}}{\Large{\wedge}} Q$  }

Let us recall the notation of minus partial order of $\mathcal{B(H)}.$ For $A,B\in \mathcal{B(H)},$ we say
$A\leq^-B$ (the symbol $\leq^-$ stands for the minus order) if there exist $P, Q\in \mathcal{B(H)}^{Id} $ such that $A = PB$ and
$A^* = QB^*.$ It follows from the above definition that  $A\leq^-B$ if and only if $A^*\leq^-B^*.$ Furthermore, [5, Proposition 3.2] implies that $A\leq^-B$ if and only if there exist $P'\in \mathcal{B(H)}^{Id} $ such that
$A = P'B$ and $R(A)\subseteq R(B).$ Thus the minus order $A\leq^-B$ induces the inclusions $R(A)\subseteq R(B)$ and $N(B)\subseteq N(A).$ The following lemma shows that the other direction holds for $A,B\in \mathcal{B(H)}^{Id}.$

{\bf Lemma 2.1.} Let $P,Q\in \mathcal{B(H)}^{Id}.$ Then the following statements are equivalent:

$(i)$   $P\preceq Q;$

$(ii)$  $R(P)\subseteq R(Q)$ and $N(Q)\subseteq N(P);$

$(iii)$   $P^{*}\preceq Q^{*};$

$(iv)$   $(I-Q) \preceq (I-P).$

{\bf Proof.} $(i)\Rightarrow (ii)$ and $(iii)\Longleftrightarrow (i)$ are straightforward.

$(ii)\Rightarrow (i).$ Since $R(P)\subseteq R(Q),$
then $QP=P.$
On the other hand,  $N(Q)\subseteq N(P)$ implies $R(P^{*})\subseteq R(Q^{*}),$ so $Q^{*}P^{*}=P^{*},$  which yields $PQ=P.$
Then $P\preceq Q.$

$(i)\Longleftrightarrow (iv).$ $(I-Q)\preceq(I-P)$ if and only if $(I-Q)(I-P)=(I-P)(I-Q)=I-Q,$ and this is the case if and only if $PQ=QP=P,$ or equivalently, $P\preceq{Q}.$
$\Box$

The following lemma is obvious from the definition.

{\bf Lemma 2.2.} Let $P\in \mathcal{P(H)}$ and $Q\in \mathcal{B(H)}^{Id}.$
Then

$(a)$ $P\preceq Q\Longleftrightarrow P\preceq Q_{or}\Longleftrightarrow P\leq Q_{or}.$

$(b)$ $Q\preceq P\Longleftrightarrow Q^{or}\preceq P \Longleftrightarrow Q^{or}\leq P.$

The following proposition give the equivalence between the existence of $P\underset{\tiny{\preceq}}{\Large{\vee}} Q$
and the existence of  $(I-P)\underset{\tiny{\preceq}}{\Large{\wedge}} (I-Q).$

{\bf Proposition 2.3.} Let $P,Q\in \mathcal{B(H)}^{Id}.$ Then the following statements are equivalent:

$(i)$   $P\underset{\tiny{\preceq}}{\Large{\vee}} Q=Q_{0};$

$(ii)$  $P^{*}\underset{\tiny{\preceq}}{\Large{\vee}} Q^{*}=Q_{0}^{*};$

$(iii)$   $(I-P)\underset{\tiny{\preceq}}{\Large{\wedge}} (I-Q)=I-Q_{0}.$

{\bf Proof.} $(i)\Rightarrow (ii).$
If $P\underset{\tiny{\preceq}}{\Large{\vee}} Q=Q_{0},$
then $P\preceq Q_{0}$ and $Q\preceq Q_{0},$
which imply $P^{*}\preceq Q_{0}^{*}$ and $Q^{*}\preceq Q_{0}^{*}$ from Lemma 2.1.

Let $Q^{'}\in \mathcal{B(H)}^{Id}$ satisfy that $P^{*}\preceq Q^{'}$ and $Q^{*}\preceq Q^{'}.$
By Lemma 2.1, we get
$$P\preceq (Q^{'})^{*} \hbox{   }\hbox{ and }\hbox{   } Q\preceq (Q^{'})^{*},$$
so $Q_{0}=P\underset{\tiny{\preceq}}{\Large{\vee}} Q\preceq (Q^{'})^{*}.$
Using Lemma 2.1 again, we conclude that $Q_{0}^{*}\preceq Q^{'},$
which induces $P^{*}\underset{\tiny{\preceq}}{\Large{\vee}} Q^{*}$ exists and $P^{*}\underset{\tiny{\preceq}}{\Large{\vee}} Q^{*}=Q_{0}^{*}.$

$(ii)\Rightarrow (iii).$
If $P^{*}\underset{\tiny{\preceq}}{\Large{\vee}} Q^{*}=Q_{0}^{*},$
then $P^{*}\preceq Q_{0}^{*}$ and $Q^{*}\preceq Q_{0}^{*}.$
Thus Lemma 2.1 implies
$$(I-Q_{0})\preceq (I-P) \hbox{   }\hbox{ and }\hbox{   }(I-Q_{0})\preceq (I-Q).$$

On the other hand, if $Q^{'}\in \mathcal{B(H)}^{Id}$ satisfy that $Q^{'}\preceq (I-P)$ and $Q^{'}\preceq (I-Q),$
then
$$P^{*}\preceq (I-Q^{'})^{*} \hbox{   }\hbox{ and }\hbox{   } Q^{*}\preceq (I-Q^{'})^{*}$$ follows from Lemma 2.1.
Thus $Q_{0}^{*}=P^{*}\underset{\tiny{\preceq}}{\Large{\vee}} Q^{*}\preceq (I-Q^{'})^{*},$
which yields $Q^{'}\preceq (I-Q_{0}).$
Thus  $(I-P)\underset{\tiny{\preceq}}{\Large{\wedge}} (I-Q)$ exists with  $(I-P)\underset{\tiny{\preceq}}{\Large{\wedge}} (I-Q)=I-Q_{0}.$
$(iii)\Rightarrow (i)$ follows in a similar way as $(ii)\Rightarrow (iii).$
 $\Box$

{\bf Lemma 2.4.}
 Let $\mathcal{M}\subseteq\mathcal{H}$ be finite dimensional and $\mathcal{N}\subseteq\mathcal{H}$ be a closed subspace.

 $(a)$ $\mathcal{M}+\mathcal{N}$ is
 a closed subspace.

 $(b)$ If $\mathcal{M}\cap\mathcal{N}=\{0\},$
then $\mathcal{M}\dotplus [\mathcal{N}\oplus(\mathcal{M}^\perp\cap \mathcal{N}^\perp)]=\mathcal{H}.$

{\bf Proof.} $(a)$ follows from [4].

$(b)$  Setting $\mathcal{M}'=\mathcal{M}^\perp\cap \mathcal{N}^\perp,$
we get from $(a)$ that $\mathcal{M}+\mathcal{N}=\overline{\mathcal{M}+\mathcal{N}}=\mathcal{M}'^\perp,$ so $\mathcal{M}+(\mathcal{N}\oplus\mathcal{M}')=\mathcal{H}.$
We claim that $\mathcal{M}\cap (\mathcal{N}\oplus \mathcal{M}^{'})=\{0\}.$
Indeed, suppose that $x\in \mathcal{M}$ and $x\in \mathcal{N}\oplus \mathcal{M}^{'}.$
Then $x=y+z,$ where $y\in \mathcal{N}$ and $z\in \mathcal{M}^{'},$
  so $z=x-y\in \mathcal{M}+\mathcal{N},$ which yields $z=0.$
Thus $x=y\in \mathcal{N}.$ Then $x=0$ follows from $\mathcal{M}\cap\mathcal{N}=\{0\}.$
Hence $\mathcal{M}\dotplus [\mathcal{N}\oplus(\mathcal{M}^\perp\cap \mathcal{N}^\perp)]=\mathcal{H}.$
$\Box$

 The following theorem 2.5 and 2.9 give an equivalent condition
 for the existences of $P \underset{\tiny{\preceq}}{\Large{\vee}} Q$ and $P \underset{\tiny{\preceq}}{\Large{\wedge}} Q,$ respectively.
 In the finite dimensional case, the existence of $P \underset{\tiny{\preceq}}{\Large{\vee}} Q$  has been
considered in [18, Lemma 3.1]. We shall extend the result to the infinite dimensional
Hilbert space.

{\bf Theorem 2.5.} Let $P,Q\in \mathcal{B(H)}^{Id}. $ Then

$(i)$ $P \underset{\tiny{\preceq}}{\Large{\vee}} Q=I$ if and only if $N(P)\cap N(Q)\subseteq \overline{R(P)+R(Q)}.$

$(ii)$ $P \underset{\tiny{\preceq}}{\Large{\wedge}} Q=0$ if and only if $R(P)\cap R(Q)\subseteq \overline{N(P)+N(Q)}.$

{\bf Proof.} $(i)$ Sufficiency. Let $Q'\in \mathcal{B(H)}^{Id}$ satisfy $Q\preceq Q^{'}$ and $P\preceq Q^{'}.$
Using Lemma 2.1, we have
$$\overline{R(P)+R(Q)}\subseteq R(Q^{'}) \hbox{   }\hbox{ and }\hbox{   }N(Q^{'})\subseteq N(P)\cap N(Q).$$
If $N(P)\cap N(Q)\subseteq \overline{R(P)+ R(Q)},$ then $N(Q^{'})\subseteq R(Q^{'}),$
and so $N(Q^{'})=\{0\}.$ Thus $Q^{'}=I,$   which yields $P \underset{\tiny{\preceq}}{\Large{\vee}} Q=I.$

Necessity. Let us assume the opposite $N(P)\cap N(Q)\nsubseteq \overline{R(P)+ R(Q)},$ and see what happens.
Then there exists
$0\neq x\in N(P)\cap N(Q)$ and $ x\notin \overline{R(P)+ R(Q)}.$
Setting  $\mathcal{M}=\overline{\{x\}},$
 we conclude that $\mathcal{M}\cap \overline{R(P)+ R(Q)}=\{0\}.$
Let $\mathcal{M}^{'}=(\mathcal{M}+\overline{R(P)+R(Q)})^\perp.$ Then Lemma 2.4 implies
$$\mathcal{M}\dotplus (\overline{R(P)+ R(Q)}\oplus \mathcal{M}^{'})=\mathcal{H}.$$
Let $Q_{0}$  be the idempotent with $R(Q_0)=\overline{R(P)+ R(Q)}\oplus {\mathcal{M}^{'}}$ and $N(Q_0)=\mathcal{M},$
then we get from Lemma 2.1 that $P\preceq Q_{0}$ and  $Q\preceq Q_{0}.$
 However $Q_{0}\neq I$ follows from $N(Q_{0})\neq 0.$
It is a contradiction with the assumption  $P\underset{\tiny{\preceq}}{\Large{\vee}} Q=I.$

$(ii)$ follows from above $(i)$ and Proposition 2.3.
$\Box$

The following proposition presents the relation between the existence of $P\underset{\tiny{\preceq}}{\Large{\vee}} Q$ and $P^{or}\underset{\tiny{\preceq}}{\Large{\vee}} Q^{or}.$

{\bf Proposition 2.6.} Let $P,\ Q\in \mathcal{B(H)}^{Id}. $ Then

$(a)$ $P^{or}\underset{\tiny{\preceq}}{\Large{\vee}} Q^{or}$ exists and $P^{or}\underset{\tiny{\preceq}}{\Large{\vee}} Q^{or}=P^{or}\vee Q^{or},$
where $P^{or}\vee Q^{or}$ is the orthogonal projection onto
the closed subspace $\overline{R(P^{or})+ R(Q^{or})}.$

$(b)$ $P\underset{\tiny{\preceq}}{\Large{\vee}} Q$ exists with $P\underset{\tiny{\preceq}}{\Large{\vee}} Q\in \mathcal{P(H)}$ if and only if $P\underset{\tiny{\preceq}}{\Large{\vee}} Q=P^{or}\underset{\tiny{\preceq}}{\Large{\vee}} Q^{or}.$

{\bf Proof.} $(a)$ Obviously, $P^{or}\preceq P^{or}\vee Q^{or}$ and $Q^{or}\preceq P^{or}\vee Q^{or}.$
Let $Q'\in \mathcal{B(H)}^{Id}$ satisfy $P^{or}\preceq Q^{'}$ and $Q^{or}\preceq Q^{'}.$
Then Lemma 2.2 implies
$P^{or}\leq Q^{'}_{or}$ and $Q^{or}\leq Q^{'}_{or},$ so
$P^{or}\vee Q^{or}\preceq Q^{'}_{or},$ which yields $P^{or}\vee Q^{or}\preceq Q^{'}.$
That is $P^{or}\underset{\tiny{\preceq}}{\Large{\vee}} Q^{or}=P^{or}\vee Q^{or}.$

$(b)$ Sufficiency is clearly.

Necessity. If $P\underset{\tiny{\preceq}}{\Large{\vee}}Q=Q_{0}\in \mathcal{P(H)},$
then $P^{or}\preceq Q_{0}$ and $Q^{or}\preceq Q_{0},$ so
  $P^{or}\vee Q^{or}\preceq Q_{0}$ follows from Lemma 2.2. On the other hand,
$P\preceq P^{or}$ and $Q\preceq Q^{or}$ imply $P, Q\preceq P^{or}\vee Q^{or}.$
Thus $P \underset{\tiny{\preceq}}{\Large{\vee}}Q=Q_{0}\preceq P^{or}\underset{\tiny{\preceq}}{\Large{\vee}} Q^{or},$ which yields $Q_{0}= P^{or}\underset{\tiny{\preceq}}{\Large{\vee}} Q^{or}.$
$\Box$

{\bf Corollary 2.7.} Let $P\in \mathcal{B(H)}^{Id}. $ Then

$(i)$ $P\underset{\tiny{\preceq}}{\Large{\vee}} (I-P)=I.$

$(ii)$  $P\underset{\tiny{\preceq}}{\Large{\vee}} P^{*}=P^{or}.$

$(iii)$ $P\underset{\tiny{\preceq}}{\Large{\wedge}} (I-P)=0.$

{\bf Proof.} $(i)$  It is straightforward that $N(P)\cap N(I-P)=N(P)\cap R(P)=0,$
which implies $N(P)\cap N(Q)\subseteq \overline{R(P)\vee R(Q)}.$
Then $P\underset{\tiny{\preceq}}{\Large{\vee}} (I-P)=I$ follows from Theorem 2.5.

$(ii)$ Clearly, $P\preceq P^{or}$ implies $P^{*}\preceq P^{or}.$
Let $Q\in \mathcal{B(H)}^{Id}$ satisfy that $P\preceq Q$ and $P^{*}\preceq Q.$
Then
$$QP=PQ=P \hbox{   }\hbox{ and }\hbox{   } P^{*}Q=QP^{*}=P^{*},$$
so $Q(P+P^{*})=(P+P^{*})Q=P+P^{*},$
which yields $R(P+P^{*})\subseteq R(Q)\cap R(Q^{*}).$
Thus $QP^{or}=P^{or}=Q^*P^{or}$ follows from following Proposition 3.3 $(iii),$
which says $P^{or}\preceq Q.$
Then $P\underset{\tiny{\preceq}}{\Large{\vee}} P^{*}=P^{or}$ as desired.

$(iii)$ follows from Theorem 2.5 $(ii).$
$\Box$

{\bf Corollary 2.8.} Let $P\in\mathcal{B(H)}^{Id}. $ Then $P_{\overline{R(P+P^{*}})}\underset{\tiny{\preceq}}{\Large{\vee}} P_{\overline{R(2I-P-P^{*}})}=I.$

{\bf Proof.} By Corollary 2.7 (i) and Proposition 2.6 (b), we know that $P^{or}\underset{\tiny{\preceq}}{\Large{\vee}}(I-P)^{or}=I.$ Then
Proposition 3.3 $(iii)$ implies $P_{\overline{R(P+P^{*}})}\underset{\tiny{\preceq}}{\Large{\vee}} P_{\overline{R(2I-P-P^{*}})}=I.$
$\Box$

The following theorem characterize an equivalent condition for the existence
of $P \underset{\tiny{\preceq}}{\Large{\vee}} Q$ with $P \underset{\tiny{\preceq}}{\Large{\vee}} Q\neq I.$

{\bf Theorem 2.9.} Let $P,Q\in \mathcal{B(H)}^{Id}.$ Then

$(i)$ $P \underset{\tiny{\preceq}}{\Large{\vee}} Q$ exists
 and $P \underset{\tiny{\preceq}}{\Large{\vee}} Q\neq I$
if and only if $N(P)\cap N(Q)\neq\{0\}$ and $(N(P)\cap N(Q))\dotplus\overline{R(P)+R(Q)}=\mathcal{H}.$
In this case, $P \underset{\tiny{\preceq}}{\Large{\vee}} Q$ is the idempotent operator onto $\overline{R(P)+R(Q)}$ along $N(P)\cap N(Q).$

$(ii)$ $P \underset{\tiny{\preceq}}{\Large{\wedge}} Q$ exists
 and $P \underset{\tiny{\preceq}}{\Large{\wedge}} Q\neq 0$
if and only if $R(P)\cap R(Q)\neq\{0\}$ and $(R(P)\cap R(Q))\dotplus\overline{N(P)+N(Q)}=\mathcal{H}.$ In this case, $P \underset{\tiny{\preceq}}{\Large{\wedge}} Q$ is the idempotent operator onto $R(P)\cap R(Q)$ along $\overline{N(P)+N(Q)}.$

{\bf Proof.} $(i)$ Sufficiency.
Let $Q_{1}\in \mathcal{B(H)}^{Id}$ with \begin{equation}N(Q_{1})= N(P)\cap N(Q) \hbox{ } \hbox{ and } \hbox{  }R(Q_{1})=\overline{R(P)+ R(Q)}.\end{equation}
Then the assumption of $N(P)\cap N(Q)\neq\{0\}$ induces $Q_1\neq I.$ Using Lemma 2.1,  we conclude from equation (2.1) that $P,Q\preceq Q_{1}.$ If
$Q'\in \mathcal{B(H)}^{Id}$ satisfies $P,Q\preceq Q',$
then \begin{equation}N(Q')\subseteq N(P)\cap N(Q)  \hbox{ } \hbox{ and } \hbox{ } R(Q')\supseteq\overline{R(P)+R(Q)}\end{equation} follow from Lemma 2.1.
Combining (2.1) and (2.2), we get that $Q_{1}\preceq Q',$ which yields
$P \underset{\tiny{\preceq}}{\Large{\vee}} Q=Q_1\neq I.$

Necessity. Let $P\underset{\tiny{\preceq}}{\Large{\vee}}  Q=Q_{2}$ and $Q_{2}\neq I.$
Then we know that $N(Q_{2})\neq\{0\},$ \begin{equation}N(Q_{2})\subseteq N(P)\cap N(Q)
\hbox{ } \hbox{ and } \hbox{ }  \overline{R(P)+ R(Q)}\subseteq R(Q_{2}).\end{equation} We claim that
 \begin{equation} N(Q_{2})= N(P)\cap N(Q)\hbox{   }\hbox{ and }\hbox{   }R(Q_{2})=\overline{R(P)+ R(Q)}.\end{equation}

Conversely, if $N(Q_{2})\neq N(P)\cap N(Q),$ then $ N(Q_{2}) \varsubsetneqq N(P)\cap N(Q),$ so there exists $y\neq 0$ such that $y\in N(P)\cap N(Q)$ and $y\notin N(Q_{2}).$

{\bf Case 1.} Suppose that $Q_{2}y\neq y.$ Then $y\notin R(Q_{2}),$ so (2.3) yields
$$\overline{\{y\}} \cap \overline{R(P)+ R(Q)}\subseteq \overline{\{y\}}\cap R(Q_{2})=\{0\}.$$
Setting $$M^{'}=(\overline{\{y\}} + \overline{R(P)+ R(Q)})^\perp,$$  we conclude from Lemma 2.4 that
$$\overline{\{y\}}\dotplus ( \overline{R(P)+ R(Q)}\oplus M^{'})=\mathcal{H}.$$
Let $Q_{3}$ be the idempotent onto $\overline {R(P)+R(Q)}\oplus {M^{'}}$ along $\overline{\{y\}}.$
Then Lemma 2.1 implies $P\preceq Q_{3}$ and $Q\preceq Q_{3}.$ However, $Q_{2}\npreceq Q_{3}$ because $N(Q_{3})=\overline{\{y\}}\nsubseteq N(Q_{2}).$
It is a contradiction with $P\underset{\tiny{\preceq}}{\Large{\vee}} Q=Q_{2}.$
Thus $N(Q_{2})=N(P)\cap N(Q)$ as desired.

{\bf Case 2.} Suppose that $Q_{2}y=y.$  As $N(Q_{2})\neq\{0\},$ we take a vector $0\neq z\in N(Q_{2}).$ Then
\begin{equation}y+z\in N(P)\cap N(Q)\hbox{   }\hbox{ and }\hbox{   } y+z\notin N(Q_{2}).\end{equation}
 Moreover, \begin{equation}Q_{2}(y+z)=Q_{2}(y)=y\neq{y+z} \end{equation}
Combining (2.5) and (2.6),
we get a contradiction by replacing $y$ with $y+z$ as in Case 1.
Thus we have $N(Q_{2})= N(P)\cap N(Q),$ which yields
\begin{equation}(N(P)\cap N(Q))\cap \overline{R(P)+ R(Q)}\subseteq N(Q_{2})\cap R(Q_{2})=\{0\}.\end{equation}

 In the following, we show that $R(Q_{2})=\overline{R(P)+ R(Q)}.$ Conversely, assume that $\overline{R(P)+ R(Q)} \varsubsetneqq R(Q_{2}).$ Then there exists $x\neq 0$ such that $x\in R(Q_2)$ and $x\notin \overline{R(P)+ R(Q)},$ so $$\overline{\{x\}}\cap\overline{R(P)+ R(Q)}=\{0\}.$$
 Setting a subspace  $$\mathcal{M}^{''}:=R(Q_2)\ominus(\overline{\{x\}}+ \overline{R(P)+ R(Q)}),$$ we know that
 $$[(\overline{\{x\}} +\overline{R(P)+ R(Q)})\oplus \mathcal{M}^{''}]\dotplus (N(P)\cap N(Q))=R(Q_{2})\dotplus N(Q_{2})=\mathcal{H},$$ which implies $$[\overline{\{x+w\}}+(\overline{R(P)+ R(Q)}\oplus\mathcal{M}^{''})]\dotplus (N(P)\cap N(Q))=\mathcal{H},$$ where $0\neq w\in N(P)\cap N(Q)=N(Q_{2}).$
Let $Q_{4}$  be the idempotent onto $\overline{\{x+w\}}+ (\overline{R(P)+ R(Q)}\oplus\mathcal{M}^{''})$
along $N(P)\cap N(Q),$ then
we conclude from Lemma 2.1 that $P\preceq Q_{4}$ and $Q\preceq Q_{4}.$
 On the other hand, it is easy to check that $x\notin R(Q_{4}),$
 so $R(Q_{2})\nsubseteq R(Q_{4}),$ which yields $Q_{2}\npreceq Q_4.$
It is a contradiction with the fact
$P\underset{\tiny{\preceq}}{\Large{\vee}}  Q=Q_{2}.$
Thus $R(Q_{2})=\overline{R(P)+ R(Q)}.$ Then (2.4) holds, which induces
$$(N(P)\cap N(Q))\dotplus\overline{R(P)+ R(Q)}=N(Q_2)\dotplus R(Q_2)=\mathcal{H}$$ and $P \underset{\tiny{\preceq}}{\Large{\vee}} Q=Q_2$ is the idempotent operator onto $\overline{R(P)+R(Q)}$ along $(N(P)\cap N(Q)).$

$(ii)$ follows from above $(i)$ and Proposition 2.3.
$\Box$

{\bf Corollary 2.10.} Let $P,\ Q\in \mathcal{B(H)}^{Id}. $
Then the following statements are equivalent:

$(i)$ $P\underset{\tiny{\preceq}}{\Large{\vee}} Q$ exists with $P\underset{\tiny{\preceq}}{\Large{\vee}} Q\in \mathcal{P(H)}\backslash\{I\};$

$(ii)$
$\overline{R(PP^{*}+ QQ^{*})}=\overline{R(P^{*}P+ Q^{*}Q)}\neq \mathcal{H};$

$(iii)$
$\{0\}\neq N(P)\cap N(Q)\subseteq N(P^{*}+P)\cap N(Q^{*}+Q)$ and
$N(P^*)\cap N(Q^*)\subseteq N(P^{*}+P)\cap N(Q^{*}+Q).$

{\bf Proof.} $(i)\Longleftrightarrow (ii).$ By the proof of Theorem 2.9, we get that $P\underset{\tiny{\preceq}}{\Large{\vee}} Q\in \mathcal{P(H)} \backslash \{I\}$
if and only if $$N(P)\cap N(Q)\neq\{0\} \hbox {  }\hbox { and }\hbox {   } (N(P)\cap N(Q))\oplus \overline{R(P)+ R(Q)}=\mathcal{H},$$
which is equivalent to    \begin{equation}0\neq N(P)\cap N(Q)=\overline{R(P)+ R(Q)}^{\perp}=N(P^{*})\cap N(Q^{*}).\end{equation}
Since $$N(P)\cap N(Q)=N(P^{*}P+ Q^{*}Q)
\hbox {  }\hbox { and }\hbox {   }  N(P^{*})\cap N(Q^{*})= N(PP^{*}+ QQ^{*}),$$
this is the case if and only if
$$
\begin{array}{rl}
\overline{R(PP^{*}+ QQ^{*})}
&=(N(P P^{*}+ QQ^{*}))^{\perp}=(N(P^{*})\cap N(Q^{*}))^{\perp}
\\&=(N(P)\cap N(Q))^{\perp}=(N(P^{*}P+Q^{*}Q))^{\perp}
\\&=\overline{R(P^{*}P+ Q^{*}Q)}\neq \mathcal{H}.
\end{array}
$$

$(i)\Longleftrightarrow (iii).$
  We observe that $N(P^{*}+P)=N(P)\cap N(P^{*}).$ Indeed, if $x\in N(P^{*}+P),$ then $$0=(P^{*}+P)^{2}x=(P^{*}+P^{*}P+PP^{*}+P)x,$$ which yields $(P^{*}P+PP^{*})x=0,$
and so $x\in N(P)\cap N(P^{*}).$ This implies that  $N(P^{*}+P)\subseteq N(P)\cap N(P^{*}).$
The other inclusion $N(P^{*}+P)\supseteq N(P)\cap N(P^{*})$ is clear.
Therefore, $$\{0\}\neq N(P)\cap N(Q)\subseteq N(P^{*}+P)\cap N(Q^{*}+Q)$$ if and only if $$\{0\}\neq N(P)\cap N(Q)\subseteq N(P^{*})\cap N(Q^{*}).$$ Similarly,
$N(P^*)\cap N(Q^*)\subseteq N(P^{*}+P)\cap N(Q^{*}+Q)$ if and only if $$ N(P^*)\cap N(Q^*)\subseteq N(P)\cap N(Q).$$ Then the assumption of (iii)
  is equivalent to (2.8), which implies $(i)\Longleftrightarrow (iii)$ as desired.
$\Box$

{\bf Corollary 2.11.} Let $P,Q\in \mathcal{B(H)}^{Id}$ and $J$ be a symmetry.
If $P$ and $Q$ are commutative with $J$ and $P\underset{\tiny{\preceq}}{\Large{\vee}} Q$ exists,
then $P\underset{\tiny{\preceq}}{\Large{\vee}} Q$ is commutative with $J$ and $P\underset{\tiny{\preceq}}{\Large{\vee}} Q=\underset{\tiny{\preceq}}{\Large{min}} \{Q': P,Q\preceq Q' \hbox{ and } Q' \hbox{ is commutative with } J\}.$

{\bf Proof.} As $J$ is a symmetry, we conclude that $J$ has the operator matrix form
$$ J=\left(\begin{array}{cc}
I&0\\0&-I
\end{array}\right):N(I-J)\oplus N(I+J),$$ so $$ P=\left(\begin{array}{cc}
P_1&0\\0&P_2
\end{array}\right):N(I-J)\oplus N(I+J), \hbox{   }\hbox{ and }\hbox{  } Q=\left(\begin{array}{cc}
Q_1&0\\0&Q_2
\end{array}\right):N(I-J)\oplus N(I+J)$$
 follows from the assumption that $P$ and $Q$ are commutative with $J,$
 where $P_1,Q_1\in\mathcal{B}(N(I-J))^{Id}$ and $P_2,Q_2\in\mathcal{B}(N(I+J))^{Id}.$
 Moreover, the existence of $P\underset{\tiny{\preceq}}{\Large{\vee}} Q,$ Theorem 2.5 and 2.9 imply that
$P_i\underset{\tiny{\preceq}}{\Large{\vee}} Q_i$ exists for $i=1,2$ and
$$P\underset{\tiny{\preceq}}{\Large{\vee}}
Q=(P_1\underset{\tiny{\preceq}}{\Large{\vee}} Q_1)\oplus (P_2\underset{\tiny{\preceq}}{\Large{\vee}} Q_2).$$
Thus $(P\underset{\tiny{\preceq}}{\Large{\vee}} Q)J=J(P\underset{\tiny{\preceq}}{\Large{\vee}} Q)$ and $$P\underset{\tiny{\preceq}}{\Large{\vee}} Q=\underset{\tiny{\preceq}}{\Large{min}} \{Q': P,Q\preceq Q' \hbox{ and } Q' \hbox{ is commutative with } J\}.$$    $\Box$

\section{ Properties of $Q_{or}$ and $Q^{or}$ }

In this section, we consider properties of the $ Q_{or}$ and $ Q^{or}.$ To show our main results, the following two lemmas are needed.

{\bf Lemma 3.1.} Let $Q\in \mathcal{B(H)}^{Id}$ and $\mathcal{M}=R(Q)\cap R(Q^{*}).$
Then $Q$ has the following operator matrix form
\begin{equation}Q=\left(\begin{array}{ccc}
I_{1}&0&0\\0&I_{2}&Q_{1}\\0&0&0\end{array}\right):\mathcal{M}\oplus (R(Q)\ominus \mathcal{M})\oplus R(Q)^{\perp},
\end{equation}
where $Q_{1}\in \mathcal{B}(R(Q)^{\perp},(R(Q)\ominus M))$ has dense range.

{\bf Proof.} It is easy to check that $\mathcal{M}$ is a reducing subspace of $Q$ and
$Q\mid_\mathcal{M}=I.$ Thus $Q$ has the operator matrix form (3.1).
If $y\in R(Q)\ominus M$ and $Q_{1}^{*}y=0,$ then $Qy=y=Q^*y$
which yields $y\in\mathcal{M},$ and hence $y=0.$
This implies that $N(Q_{1}^{*})=0,$ so $Q_{1}$ has dense range.
$\Box$

{\bf Lemma 3.2.} Let $A\in \mathcal{B(K,H)}$ and
$\widetilde{A}:=\left(\begin{array}{cc}
I&A\\A^{*}&A^{*}A\end{array}\right).$
Then $\widetilde{A}\in \mathcal{B(\mathcal{H}\oplus\mathcal{K})}^{+}$ and
\begin{equation}\widetilde{A}^{\frac{1}{2}}=\left(\begin{array}{cc}
(I+AA^{*})^{-\frac{1}{2}}&(I+AA^{*})^{-\frac{1}{2}}A\\
A^{*}(I+AA^{*})^{-\frac{1}{2}}&(I+A^{*}A)^{-\frac{1}{2}}A^{*}A\end{array}\right).\end{equation}

{\bf Proof.} It is a direct verification.
$\Box$

The following proposition gives some specific structures
of $Q_{or}$ and $Q^{or}.$

{\bf Proposition 3.3.}  Let $Q\in \mathcal{B(H)}^{Id}.$ Then

$(i)$  $Q_{or}=P_{R(Q)\cap R(Q^{*})}.$

$(ii)$ $Q^{or}=I-(I-Q)_{or}.$

$(iii)$ $Q^{or}=P_{N(Q+Q^{*})^{\perp }}.$

{\bf Proof.} $(i)$ By Lemma 3.1,  we get that $$QP_{R(Q)\cap R(Q^{*})}=P_{R(Q)\cap R(Q^{*})}=P_{R(Q)\cap R(Q^{*})}Q,$$  so $P_{R(Q)\cap R(Q^{*})}\preceq Q.$
On the other hand, if $ \mathcal{P(H)}\ni P\preceq Q,$ then $PQ=QP=P,$ which implies $R(P)\subseteq R(Q)\cap R(Q^{*}).$ Thus $P\preceq P_{R(Q)\cap R(Q^{*})},$ so $$ Q_{or}=\underset{\tiny{\preceq}}{\Large{max}} \{P: P\preceq Q \hbox{ and } P\in \mathcal{P(H)}\}=P_{R(Q)\cap R(Q^{*})}.$$

$(ii)$ is trivial from Lemma 2.1 and the definitions of $Q_{or}$ and $Q^{or}.$

$(iii)$ Using $(ii),$ we know that $$Q^{or}=I-(I-Q)_{or}=I-P_{R(I-Q)\cap R(I-Q^{*})}=P_{(N(Q)\cap N(Q^{*}))^{\perp }}=P_{N(Q+Q^{*})^{\perp }}.$$
$\Box$

{\bf Corollary 3.4.} Let $Q\in \mathcal{B(H)}^{Id}$ and $J$ be a symmetry. If $Q$ is a $J$-projection, then

$(i)$  $\underset{\tiny{\preceq}}{\Large{max}} \{P: P\preceq Q,\ P\in \mathcal{P(H)},\ P \hbox{ is } J\hbox{-projection}\}=P_{R(Q)\cap R(Q^{*})}.$

$(ii)$  $\underset{\tiny{\preceq}}{\Large{min}} \{P: Q\preceq P,\ P\in \mathcal{P(H)},\ P\hbox{ is } J\hbox{-projection}\}=P_{N(Q+Q^{*})^{\perp }}.$

{\bf Proof.} $(i)$ By Proposition 3.3 $(i),$ we only need to show that $P_{R(Q)\cap R(Q^{*})}$ is a $J$-projection, that is $JP_{R(Q)\cap R(Q^{*})}=P_{R(Q)\cap R(Q^{*})}J.$
 Let $x\in {R(Q)\cap R(Q^{*})}.$ Then
$Qx=Q^{*}x=x,$ and since $JQ=Q^{*}J,$ so we have
$$QJx=JQ^{*}x=Jx \hbox{   }\hbox{ and }\hbox{   } Q^{*}Jx=JQx=Jx.$$
Thus $J^{*}x=Jx\in {R(Q)\cap R(Q^{*})},$
 which implies that  ${R(Q)\cap R(Q^{*})}$ is a reducing subspace of $J.$
 Hence  $JP_{R(Q)\cap R(Q^{*})}=P_{R(Q)\cap R(Q^{*})}J.$

$(ii)$ Using Proposition 3.3 $(ii),$ we need to show that $JP_{N(Q+Q^{*})^{\bot}}=P_{N(Q+Q^{*})^{\bot}}J,$ which is equivalent to $JP_{N(Q+Q^{*})}=P_{N(Q+Q^{*})}J.$
Let $x\in N(Q+Q^{*}).$ Then $JQ=Q^{*}J$ yields
$$(Q+Q^{*})Jx=(JQ^{*}+JQ)x=J(Q^{*}+Q)x=0$$
so $J^{*}x=Jx\in N(Q+Q^{*}). $
Thus $ N(Q+Q^{*})$ is a reducing subspace of $J,$
which induces $P_{N(Q+Q^{*})}J=JP_{N(Q+Q^{*})}.$
$\Box$

The following result is an extension of [14, Proposition 1].

{\bf Proposition 3.5.}  Let $Q\in \mathcal{B(H)}^{Id}.$ Then
$Q_{or}=P_{N(I-|Q|)}=P_{N(2I-Q-Q^{*})}.$

{\bf Proof.}
Suppose that $Q$ has the form as (3.1). Then by Lemma 3.2 we have
\begin{align*}\begin{array}{ll} &|Q|=(Q^{*}Q)^{\frac{1}{2}}=\left(\begin{array}{ccc}
I_{1}&0&0\\0&I_{2}&Q_{1}\\0&Q_{1}^{*}&Q_{1}^{*}Q_{1}\end{array}\right)^{\frac{1}{2}}\\=&\left(\begin{array}{ccc}
I_{1}&0&0\\0&(I_{2}+Q_{1}Q_{1}^{*})^{-\frac{1}{2}}&(I_{2}+Q_{1}Q_{1}^{*})^{-\frac{1}{2}}Q_{1}\\
0&Q_{1}^{*}(I_{2}+Q_{1}Q_{1}^{*})^{-\frac{1}{2}}&(I_{3}+Q_{1}^{*}Q_{1})^{-\frac{1}{2}}Q_{1}^{*}Q_{1}\end{array}\right)
\end{array}\end{align*}
Setting
$$\widetilde{Q}:=\left(\begin{array}{cc}
(I_{2}+Q_{1}Q_{1}^{*})^{-\frac{1}{2}}&(I_{2}+Q_{1}Q_{1}^{*})^{-\frac{1}{2}}Q_{1}\\
Q_{1}^{*}(I_{2}+Q_{1}Q_{1}^{*})^{-\frac{1}{2}}&(I_{3}+Q_{1}^{*}Q_{1})^{-\frac{1}{2}}Q_{1}^{*}Q_{1}
\end{array}\right),$$
we know that $N(I-|Q|)=(R(Q)\cap R(Q^{*}))\oplus N(\widetilde{Q}-I).$ We claim that  $N(\widetilde{Q}-I)=\{0\}.$
Indeed, if $\left[\begin{array}{c}x\\y\end{array}\right]\in  (R(Q)\ominus \mathcal{M})\oplus R(Q)^{\perp}$ satisfies
$\widetilde{Q}\left[\begin{array}{c}
x\\y\end{array}\right]=\left[\begin{array}{c}
x\\y\end{array}\right],
$
 then
\begin{equation}\begin{cases}(I_{2}+Q_{1}Q_{1}^{*})^{-\frac{1}{2}}x+
(I_{2}+Q_{1}Q_{1}^{*})^{-\frac{1}{2}}Q_{1}y=x\qquad\qquad\ \ \ \ \ \ \ \ \  \
&\\Q_{1}^{*}(I_{2}+Q_{1}Q_{1}^{*})^{-\frac{1}{2}}x+(I_{3}+Q_{1}^{*}Q_{1})^{-\frac{1}{2}}Q_{1}^{*}Q_{1}y=y.\qquad \qquad\ \ \end{cases}\end{equation}
Thus
$$x+Q_{1}y=(I_{2}+Q_{1}Q_{1}^{*})^{\frac{1}{2}}x \hbox{  }\hbox{   }\hbox{ and }\hbox{ }
 Q_{1}^{*}x+Q_{1}^{*}Q_{1}y=(I_{3}+Q_{1}^{*}Q_{1})^{\frac{1}{2}}y, $$
and hence $$(I_{3}+Q_{1}^{*}Q_{1})^{\frac{1}{2}}y=Q_1^*(I_{2}+Q_{1}Q_{1}^{*})^{\frac{1}{2}}x
=(I_{3}+Q_{1}^{*}Q_{1})^{\frac{1}{2}}Q_{1}^{*}x,$$
 which means $y=Q_{1}^{*}x.$ Using the first equation of (3.3),  we have
$$(I_{2}+Q_{1}Q_{1}^{*})^{\frac{1}{2}}x=(I_{2}+Q_{1}Q_{1}^{*})^{-\frac{1}{2}}x+
(I_{2}+Q_{1}Q_{1}^{*})^{-\frac{1}{2}}Q_{1}Q_{1}^{*}x=x,$$
 which implies $Q_{1}Q_{1}^{*}x=0.$
Since $Q_{1}^{*}$ is injective, it follows $x=0,$
which yields $y=Q_{1}^{*}x=0.$
Thus $N(\widetilde{Q}-I)=\{0\},$ so $Q_{or}=P_{N(I-|Q|)}$ follows from Proposition (i).
Furthermore, Proposition 3.3 (ii) implies
$$P_{N(2I-Q-Q^{*})}= I-(I-Q)^{or}=Q_{or}.$$     $\Box$

{\bf Lemma 3.6.} Let $P,\ Q\in \mathcal{B(H)}^{Id}.$ Then $P^{or}\preceq Q_{or} $ if and only if $Q=P^{or}+Q_{1},$ where $Q_{1}\in \mathcal{B(H)}^{Id}$ and $P^{or}Q_{1}=Q_{1}P^{or}=0.$

{\bf Proof.}  Sufficiency is straightforward.

Necessity.  Let $Q_1=Q-P^{or}.$ Then $$Q_1^2=(Q-P^{or})^2=Q^2-QP^{or}-P^{or}Q+P^{or}=Q-P^{or}=Q_1$$
and $$P^{or}Q_{1}=P^{or}Q-P^{or}=0=QP^{or}-P^{or}
=Q_{1}P^{or}.$$   $\Box$

The following theorem characterize a necessary and sufficient condition under which $P^{or}\preceq Q_{or}$ for all $Q\in \mathcal{B(H)}^{Id} $ with $P\prec Q$
$(P\prec Q$ denotes $P\preceq Q$ and $P\neq Q).$

{\bf Theorem 3.7.}  Let $P\in \mathcal{B(H)}^{Id}.$ Then $P^{or}\preceq Q_{or}$ for all $Q\in \mathcal{B(H)}^{Id} $ with $P\prec Q$ if and only if  $P\in \mathcal{P(H)}$ or $dim R(P)^\perp\leq1.$

{\bf Proof.} Sufficiency. If $P\in \mathcal{P(H)}$ and $P\prec Q,$ then $P^{or}=P\preceq Q_{or}$ is obvious.
Furthermore,  it is easy to verify that $dim R(P)^\perp\leq1$ and $P\prec Q$ imply $P=I$ or $Q=I,$
so desired conclusion holds.

Necessity. Let us assume that $dim R(P)^\perp\geq2$ and $P\notin \mathcal{P(H)},$ and see what happens.
Let $P$ be as (3.1). That is $$P=\left(\begin{array}{ccc}
I_{1}&0&0\\0&I_{2}&P_{1}\\0&0&0\end{array}\right):\mathcal{M}\oplus (R(P)\ominus \mathcal{M})\oplus R(P)^{\perp},$$
 where $\mathcal{M}=R(P)\cap R(P^{*})$ and $P_{1}\in \mathcal{B}(R(P)^{\perp},(R(Q)\ominus M))$ has dense range.

{\bf Case 1.} $N(P_{1})=0.$ Let $Q\in\mathcal{B(H)}$ on the space decomposition $\mathcal{H}=\mathcal{M}\oplus (R(P)\ominus \mathcal{M})\oplus R(P)^{\perp}$ have the operator matrix form
$$
Q=\left(\begin{array}{ccc}
I_{1}&0&0\\0&I_{2}&P_{1}-P_{1}Q_{2}\\0&0&Q_{2}
\end{array}\right),
$$
where $Q_{2}\neq 0,I$ and $Q_{2}\in \mathcal{B}(R(P)^{\perp})^{Id}$ ($Q_{2}$ exists, as $dim R(P)^\perp\geq2).$
By a direct calculation, we get
$$Q^{2}=Q \hbox{   }\hbox{ and }\hbox{   }PQ=QP=P,$$
and hence $P\prec Q.$

On the other hand, $N(P+P^*)=0$ follows from
$N(P_{1})=0$ and $N(P_{1}^*)=0.$
And by Proposition 3.3, we have
$P^{or}=I.$
However, $(Q-P^{or})P^{or}=Q-I\neq0,$ so Lemma 3.6 yields that $P^{or}\npreceq Q_{or}.$
This is a contradiction. Hence $P \in \mathcal{P(H)}.$

{\bf Case 2.} $N(P_{1})\neq 0.$  Then $P$ on the space decomposition $\mathcal{H}=\mathcal{M}\oplus (R(P)\ominus \mathcal{M}) \oplus N(P_{1})^{\perp} \oplus N(P_{1})$ has the operator matrix form
$$
P=\left(\begin{array}{cccc}
I_{1}&0&0&0\\0&I_{2}&P_{11}&0\\0&0&0&0\\0&0&0&0
\end{array}\right),
$$
where $P_{11}\in \mathcal{B}(N(P_{1})^{\perp},R(P)\ominus \mathcal{M})$ is injective and has dense range, as $P_{1}$ has dense range.
Define $Q'\in \mathcal{B(H)}$ on the space decomposition $\mathcal{H}=\mathcal{M}\oplus (R(P)\ominus \mathcal{M}) \oplus N(P_{1})^{\perp} \oplus N(P_{1})$ by the operator matrix form
$$
Q'=\left(\begin{array}{cccc}
I_{1}&0&0&0\\0&I_{2}&0&-P_{11}Q_{11}\\0&0&I_{3}&Q_{11}
\\0&0&0&0\end{array}\right),
$$
where $0\neq Q_{11}\in \mathcal{B}(N(P_{1}), N(P_{1})^{\perp}).$
A direct calculation implies
$$Q'^{2}=Q' \hbox{   }\hbox{ and }\hbox{   }PQ'=Q'P=P,$$
which yields $P\prec Q'.$

Using Proposition 3.3 again,
we get that  $P^{or}=diag(I_1,I_2,I_3, 0),$ which yields
$P^{or}(Q'-P^{or})\neq0.$ Then Lemma 3.6 implies that $P^{or}\npreceq Q_{or}.$
This is a contradiction.
Hence, if $dim(R(P))^{\perp}\geq{2},$ then $P\in \mathcal{P(H)},$ so Necessity holds.
$\Box$

Corollary 3.4 above shows that if $Q\in \mathcal{B(H)}^{Id}$ is a $J$-projection,
then $Q^{or}\in \mathcal{P(H)}$ is a $J$-projection.
A natural problem is that
 whether there is a $J$-projection $Q\in \mathcal{B(H)}^{Id} \backslash \mathcal\{{\mathcal{P(H)}}\}$
such that $Q^{or}=P,$ if $P\in \mathcal{P(H)}$ is a $J$-projection. The following result gives the answer of this problem.

{\bf Theorem 3.8.} Let $P\in \mathcal{P(H)}$ and $J$ be a symmetry with $JP=PJ.$
 
$(i)$ There exists an idempotent $Q\in \mathcal{B(H)}^{Id} \backslash \mathcal\{{\mathcal{P(H)}}\}$
such that $Q^{or}=P$ and $JQ=Q^{*}J$ if and only if $dim R(P)\geq 2$ and $(I\pm J)P\neq 0.$

$(ii)$ There exists an idempotent $Q'\in \mathcal{B(H)}^{Id} \backslash \mathcal\{{\mathcal{P(H)}}\}$
such that $Q'_{or}=P$ and $JQ'=Q'^{*}J$ if and only if $dim R(I-P)\geq 2$ and $(I\pm J)(I-P)\neq 0.$

{\bf Proof.} $(i)$ Sufficiency.
Since $(I\pm J)P\neq 0$ and $JP=PJ,$ $J$ has the operator matrix form
$$ J=\left(\begin{array}{cc}
J_{1}&0\\0&J_{2}
\end{array}\right):R(P)\oplus R(P)^{\perp},$$
where $J_{1}$, $J_{2}$ are symmetries with $J_{1}\neq \pm I_{1}.$  Thus there exist unit vectors $x_{1},x_{2}\in R(P)$ such that $x_{1}\perp x_{2},$
$$Jx_{1}=J_{1}x_{1}=x_{1}\hbox{   }\hbox{ and }\hbox{   } Jx_{2}=J_{1}x_{2}=-x_{2},$$
so $J$ on the space decomposition $\mathcal{H}=\overline{\{x_{1}\}} \oplus \overline{\{x_{2}\}} \oplus (R(P)\ominus \overline{\{x_{1},x_{2}\}} )\oplus R(P)^{\perp}$ has the operator matrix form
$$ J=\left(\begin{array}{cccc}
1&0&0&0\\0&-1&0&0\\0&0&J_{11}&0\\0&0&0&J_{2}
\end{array}\right),$$
where $J_{11}$ is a symmetry.
 Let $Q\in \mathcal{B(H)}$ on the space decomposition $\mathcal{H}=\overline{\{x_{1}\}} \oplus \overline{\{x_{2}\}} \oplus (R(P)\ominus \overline{\{x_{1},x_{2}\}} )\oplus R(P)^{\perp}$ have the operator matrix form
$$ Q=\left(\begin{array}{cccc}
\frac{3}{2}&\frac{\sqrt{-3}}{2}&0&0\\\frac{\sqrt{-3}}{2}&-\frac{1}{2}&0&0\\0&0&I&0\\0&0&0&0
\end{array}\right).$$
Then it is easy to check that $JQ=Q^{*}J ,$ $Q\in \mathcal{B(H)}^{Id} \backslash \mathcal\{{\mathcal{P(H)}}\}$ and $N(Q+Q^{*})=R(P)^{\perp},$
so $Q^{or}=P$ follows from Proposition 3.3 (iii).

Necessity.
Suppose that $Q^{or}=P,$  which implies $QP=PQ=Q,$ and hence $R(Q)\subseteq R(P).$
If $dim R(P)=1,$ then $dim R(Q)=1.$ So $P=x\otimes x$ for a unit vector $x\in \mathcal{H},$
and $Q=y\otimes z$ for non-zero vectors $y,z\in \mathcal{H}.$
Therefore,
$$ QP=(y\otimes z)(x\otimes x)=\langle x,z\rangle(y\otimes x)=y\otimes z=Q$$
and
$$ PQ=(x\otimes x)(y\otimes z)=\langle y,x \rangle(x\otimes z)=y\otimes z=Q.$$
Thus
$$z=\langle z,x \rangle x \hbox{  } \hbox{  and } \hbox{  } \ y=\langle y,x \rangle x,$$
so $Q=y\otimes z=\lambda (x\otimes x),$  for $0\neq\lambda \in \mathbb{C}.$
Moreover, $\lambda^{2} (x\otimes x)=\lambda (x\otimes x)$ follows from $Q^{2}=Q,$
which implies  $\lambda =1.$ Hence $Q=x\otimes x\in \mathcal{P(H)},$
which is a contradiction with the fact $Q\in \mathcal{B(H)}^{Id} \backslash \mathcal\{{\mathcal{P(H)}}\},$
so $dim R(P)\geq2.$

Conversely, we assume that $(I-J)P=0.$ Then $P=JP,$  so $J$ has the operator matrix form
$$J=\left(\begin{array}{cc}I&0\\0&J'\end{array}\right):R(P)\oplus R(P)^{\perp},$$
where $J'\in {\mathcal{B}}(R(P)^{\perp})$ is a symmetry.
Let
$$ Q=\left(\begin{array}{cc}Q_{11}&Q_{12}\\Q_{21}&Q_{22}
\end{array}\right):R(P)\oplus R(P)^{\perp}.$$
 Owing to Proposition 3.3 and the equation $Q^{or}=P,$ we get that $N(Q+Q^{*})=R(P)^{\perp},$ so for all $x\in R(P)^{\perp},$ we have
$$
\begin{array}{ll}(Q+Q^{*})
\left[\begin{array}{c}0\\x\end{array}\right]
=\left(\begin{array}{cc}
Q_{11}+Q_{11}^{*}&Q_{12}+Q_{21}^{*}\\Q_{21}
+Q_{12}^{*}&Q_{22}+Q_{22}^{*}\end{array}\right)
\left[\begin{array}{c}
0\\x\end{array}\right]=0,
\end{array}
$$
 which yields
\begin{equation} Q_{12}+Q_{21}^{*}=0 \ \ \ \hbox{ and }
\ \ \ Q_{22}+Q_{22}^{*}=0. \end{equation}

On the other hand, it follows from the equation $JQ=Q^{*}J$ that
$$\left(\begin{array}{cc}
Q_{11}&Q_{12}\\J'Q_{21}&J'Q_{22}\end{array}\right)=
\left(\begin{array}{cc}
Q_{11}^{*}&Q_{21}^{*}J'\\Q_{12}^{*}&Q_{22}^{*}J'\end{array}\right),$$
 which implies
\begin{equation} Q_{11}=Q_{11}^{*}\qquad  \hbox{ and } \qquad\
 Q_{12}=Q_{21}^{*}J'. \end{equation}
Combinbing (3.4) and (3.5),  we have
$$ Q=\left(\begin{array}{cc}
Q_{11}&Q_{21}^{*}J'\\-J'Q_{21}&Q_{22}
\end{array}\right)=\left(\begin{array}{cc}
Q_{11}^*&Q_{21}^{*}J'\\-J'Q_{21}&-Q_{22}^*
\end{array}\right)\in \mathcal{B(H)}^{Id}.$$
Moreover, $Q^{2}=Q$ yields
$$ \left(\begin{array}{cc}
Q_{11}^{2}-Q_{21}^{*}Q_{21}&Q _{11}Q_{21}^{*}J'+Q_{21}^{*}J'Q_{22}
\\-J'Q_{21}Q_{11}-Q_{22}J'Q_{21}&Q_{22}^{2}-J'Q_{21}Q_{21}^{*}J'
\end{array}\right)=\left(\begin{array}{cc}
Q_{11}&Q_{12}\\Q_{21}&Q_{22}
\end{array}\right),$$
which implies
\begin{equation}
Q_{11}^{2}-Q_{21}^{*}Q_{21}=Q_{11} \qquad  \hbox{ and } \qquad\
 Q_{22}^{2}-J'Q_{21}Q_{21}^{*}J'=Q_{22}. \end{equation}

 Hence
\begin{equation}Q_{22}=Q_{22}^{2}-J'Q_{21}Q_{21}^{*}J'
=(Q_{22}^{*})^{2}-J'Q_{21}Q_{21}^{*}J'=Q_{22}^{*}.\end{equation}
Then (3.4) induces
  $Q_{22}=0$ and $J'Q_{21}Q_{21}^{*}J'=0,$
 so $Q_{21}=0.$
Thus $Q_{11}^{2}=Q_{11}$ by (3.6).
Using (3.5) again, we get that $Q_{12}=0$ and
$Q_{11}\in\mathcal{P}(R(P)),$
which means
$$Q=\left(\begin{array}{cc}
Q_{11}&0\\0&0
\end{array}\right)\in \mathcal{P(H)}.$$
This is a contradiction with the assumption $Q\in \mathcal{B(H)}^{Id} \backslash \mathcal\{{\mathcal{P(H)}}\}.$
Therefore, $(I-J)P\neq 0 $ as desired.
In a similar way, we have $(I+J)P\neq 0.$

$(ii)$ follows from above $(i)$ and Proposition 3.3 $(ii).$
$\Box$

The following result shows the specificity
of $Q-P\in\mathcal{B(H)}^{+},$ when $P\preceq Q$ for $P,\ Q\in \mathcal{B(H)}^{Id}.$

{\bf Proposition 3.9.} Let $P,\ Q\in \mathcal{B(H)}^{Id}.$ If $P\preceq Q,$ then the following statement are equivalent:

$(i)$   $Q-P\geq 0;$

$(ii)$   $Q-P$ is self-adjoint;

$(iii)$  $Q-P$ is an orthogonal projection;

$(iv)$  $Q+Q^{*}\geq P+P^{*}.$

{\bf Proof.}  $(i)\Rightarrow (ii)$ is obvious.

$(ii)\Rightarrow (iii).$ As $ PQ=QP=P,$ we know that
$$(Q-P)^{2}=(Q-P)(Q-P)=Q^{2}-QP-PQ+P^{2}=Q-P.$$
 Thus $(ii)$ implies that $Q-P$ is an orthogonal projection as desired.

$(iii)\Rightarrow (iv).$ It is clear that
$$Q+Q^{*}-(P+P^{*})=(Q-P)+(Q-P)^{*}=2(Q-P)\geq 0,$$
so $Q+Q^{*}\geq P+P^{*}.$

$(iv)\Rightarrow(i).$ Setting $A:=Q-P,$
we conclude from $(iv)$ that $A+A^{*}\geq 0.$ As $A^{2}=A,$ then $A$ as an operator on the space decomposition $R(A)\oplus R(A)^{\perp}$ has the operator matrix form
$$
A=\left(\begin{array}{cc}
I&A_{1}\\0&0\end{array}\right),
$$
which yields
$$
A+A^{*}=\left(\begin{array}{cc}
2I&A_{1}\\A_{1}^{*}&0\end{array}\right).
$$
So $A_{1}=0$ by $A+A^{*}\geq 0.$ Thus
$Q-P=A\geq 0.$
$\Box$

At last, we present a result about the continuity of the map: $P\rightarrow P^{or}.$

{\bf Proposition 3.10.} Let $Q_{n},Q\in \mathcal{B(H)}^{Id}$ and $J$ be a symmetry. Suppose that ${Q_{n}}$ is a sequence of $J$-projections. Then

$(i)$ If $Q_{n}{\overset{\tiny{WOT}}{\Large{\nearrow}}}Q,$ then $Q$ is $J$-projection and $Q_{n}^{or}{\overset{\tiny{WOT}}{\Large{\nearrow}}} Q^{or}.$

$(ii)$ If $Q_{n}{\overset{\tiny{WOT}}{\Large{\searrow}}} Q,$ then $Q$ is $J$-projection and $(Q_{n})_{or}{\overset{\tiny{WOT}}{\Large{\searrow}}} Q_{or}.$

{\bf Proof.} $(i).$  For all vectors $ x,y\in \mathcal{H},$   we have
$$\langle  JQ_{n}x,y \rangle=\langle  Q_{n}x,Jy \rangle \xrightarrow[n\rightarrow \infty] \  \langle  Qx,Jy \rangle$$
and
$$\langle  Q_{n}^{*}Jx,y \rangle=\langle Jx,Q_{n}y \rangle\ \xrightarrow[n\rightarrow \infty] \   \langle  Jx,Qy \rangle.$$
Thus $\langle Qx,Jy \rangle=\langle Jx,Qy \rangle$ follows from the fact that $Q_{n}$ are $J$-projection for $n=1,2,\cdots.$
 Then $JQ=Q^{*}J,$
and hence $Q$ is $J$-projection.

For any $n_{0}\in \mathbb{Z}^{+},$ if $n\geq n_{0},$
then $Q_{n_{0}}\preceq Q_{n}$ implies
$Q_{n_{0}}Q_{n}=Q_{n}Q_{n_{0}}=Q_{n_{0}}.$
 Thus
$$\langle  Q_{n_{0}}x,y \rangle=\langle  Q_{n_{0}}Q_{n}x,y \rangle\ \xrightarrow[n\rightarrow \infty] \   \langle Q_{n_{0}}Qx,y \rangle,$$
which implies that $\langle  Q_{n_{0}}x,y \rangle=\langle Q_{n_{0}}Qx,y \rangle.$
Analogously, we have $\langle Q_{n_{0}}x,y \rangle=\langle  QQ_{n_{0}}x,y \rangle.$
Thus $QQ_{n_{0}}=Q_{n_{0}}Q=Q_{n_{0}},$ that is $Q_{n_{0}}\preceq Q,$ and hence $Q_{n_{0}}^{or}\preceq Q^{or}.$
As $\{Q_{n}^{or}\}$ is a increasing sequence,
 then there exists an orthogonal projection $P$ such that $Q_{n}^{or}{\overset{\tiny{WOT}}{\Large{\nearrow}}}P,$
which implies $P\preceq Q^{or}.$

On the other hand,   it is clear that
$$
\begin{array}{rl}
\langle (PQ-Q)x,y \rangle
&=\langle (PQ-PQ_{n})x,y \rangle+\langle (PQ_{n}-Q)x,y \rangle
\\&=\langle  P(Q-Q_{n})x,y \rangle+\langle (PQ_{n}^{or}Q_{n}-Q)x,y \rangle
\\&=\langle  P(Q-Q_{n})x,y \rangle+\langle (Q_{n}-Q)x,y \rangle\xrightarrow[n\rightarrow \infty] \   0,
\end{array}
$$
so $PQ=Q.$
Similarly, we get that $QP=Q.$
Therefore, $Q\preceq P,$ which yields $Q^{or}\preceq P.$ Thus $P=Q^{or},$ so $Q_{n}^{or}{\overset{\tiny{WOT}}{\Large{\nearrow}}} Q^{or}.$ In a similarly way, we might show that $(ii)$ holds.
$\Box$

\end{document}